\documentclass[a4paper,11pt]{article}

\usepackage{amssymb,amsfonts,amsmath,amsfonts,amsthm}
\usepackage{graphicx}
\usepackage{hyperref}
\usepackage{mathtools}
\hypersetup{
    colorlinks = true,
    urlcolor   = blue,
    citecolor  = black,
}
\usepackage{bookmark}
\usepackage{breakurl}
\usepackage{booktabs}
\usepackage{colortbl}
\usepackage[boxed]{algorithm2e}
\usepackage{calc}
\usepackage{soul}
\usepackage{bm}

\usepackage{newtxtext}
\usepackage{newtxmath}
\usepackage{epstopdf}

\usepackage{color}
\definecolor{mycol}{cmyk}{0,1,0,0}

\usepackage[authoryear]{natbib}
\newsavebox{\astrutbox}
\sbox{\astrutbox}{\rule[-5pt]{0pt}{20pt}}

\newcommand\thalf{\ensuremath{{\textstyle\frac{1}{2}}}}

\newcommand\Rey{\mbox{\textit{Re}}}  
\providecommand{\keywords}[1]{\textit{Keywords:} #1}

\title{Energy Transfer Dynamics Generated by Non-Axisymmetric Tornado-Type Flows}
\author{
Afifah Maya Iknaningrum$^{1}$,
Pen-Yuan Hsu$^{2}$,
Tsuyoshi Yoneda$^{3}$
and
Hirofumi Notsu$^{4}$
\bigskip\\
\normalsize
$^1$
Graduate School of Natural Science and Technology, Kanazawa University, Kanazawa 920-1192, Japan \\
\normalsize
$^2$ 
Graduate School of Mathematical Science, The University of Tokyo, Tokyo 153-8914, Japan \\
\normalsize
$^3$ 
Graduate School of Economics, Hitotsubashi University, Tokyo 186-8601, Japan \\
\normalsize
$^4$
Faculty of Mathematics and Physics, Kanazawa University, Kanazawa 920-1192, Japan 
}
\date{}
\begin{document}
\maketitle
\begin{abstract}
The energy cascade in turbulence, first statistically described by Richardson (1922) and Kolmogorov (1941), lacked connection to the underlying fluid dynamics. Recent numerical studies of Goto et al. (2017) and  Yoneda et al. (2022) revealed scale-local energy transfer via vortex stretching but remained within spatial statistics. This study aims to uncover the time-dependent elementary process behind the energy cascade by constructing a tornado-type flow in a non-axisymmetric curved cylindrical domain. Our approach reveals specific vortex dynamics responsible for energy transfer, offering new insight into the physical mechanisms of turbulence.
\smallskip\\
\keywords{Navier--Stokes equations, Rotating turbulence, Vortex dynamics}
\end{abstract}

\section{Introduction}

\subsection{Energy transfer in turbulence}

The energy cascade between the hierarchy of scales was described first by \citet{Richardson1922book}, who provided a statistical description of this cascade for turbulent flows. This idea further leads to the concept of self-similar flow structures, predicted by \citet{Kolmogorov1941}, more precisely, he proposed the similarity hypothesis to explain the -5/3 power law and energy cascade, but it was based on dimensional analysis (i.e. statistical description) that was not directly related to the mathematical structure of the solutions derived from the Navier--Stokes equations. In short, he could not reach an insight into what kind of fluid behavior generates energy cascade.

With the aid of direct numerical simulation, \citet{GotoEtAl17} found that the 3D Navier--Stokes turbulence in a periodic box has a hierarchical structure of vortex tubes. At each of these scales, the vortex tubes generate a strain tensor that contributes to the generation of smaller-scale vortex tubes. After their celebrated work, \citet{YonedaEtA2022} successfully reformulated Kolmogorov's -5/3 power law (see also \citet{TsuruhashiEtAl22}). More specifically, they defined an energy transfer function in terms of vortex stretching/compressing \citep[Eq.~(2.6)]{YonedaEtA2022}, and by using it, they were able to identify a clear scale-local energy transfer structure from turbulence snapshots in spatial averages at each time \citep[Figure~1]{YonedaEtA2022}. From this result, interaction between two adjacent scales (ratio of approximately $1.7$) with vortex stretching and compressing could be expected in turbulence.

However, even in this numerical result, they have used spatial statistics and have not yet been able to clarify the elementary process as a more specific time evolution. The purpose of this study is to clarify this issue, and for the first step, we proceeded with a tornado-type flow as a typical model. More precisely, in this paper we further explore by constructing the tornado-type flow in a non-axisymmetric curved cylindrical domain that induces energy transfer to adjacent like smaller scale vortices. None of the numerous works to date attempted to find such specific vortex dynamics, so, this study will be a cornerstone for understanding the vortical structures of the energy transfer in turbulence.

The next subsection provides a more detailed explanation of the setting, history and other motivations behind this tornado-type flow.

\subsection{Tornado-type flows}

Since tornadoes represent one of the most prominent examples of vortex-dominated flows, we now turn to a brief review of observational and experimental studies of tornado dynamics. Research into tornado dynamics has employed observational (field-based), experimental (laboratory-based), and computational (numerical simulation) approaches. In observational studies, the pioneering work of \citet{FujitaEtAl70} laid the foundation for quantifying tornado intensity, which was later expanded through radar observations and storm chaser reports \citep{Bluestein13,BluesteinEtAl18}. These studies have revealed essential features of tornado structure and evolution. In parallel, laboratory-based experiments have been conducted to recreate tornado-like flows under controlled conditions \citep{Ward72,ChurchEtAl79,TariEtAl10}.

Numerical simulations complement these methods by offering cost-effective and flexible tools for investigating tornado behavior, provided that initial and boundary conditions are appropriately defined. Many simulations are conducted in axisymmetric domains \citep{Rotunno77,HsuNotsuYoneda16} and typically evaluate the sensitivity to flow parameters such as swirl ratio \citep{NolanFarrell98,LiuEtAl20,ZhaoEtAl23}. Other works aim to replicate laboratory experiments in numerical frameworks \citep{YuanEtAl19} to investigate swirl ratio \citep{IshiharaEtAl11}, near-ground behavior \citep{KuaiEtAl08}, surface roughness, and transition phenomena \citep{LiuIshihara16}. Additionally, several theoretical and analytical studies have been conducted to provide a broader understanding of tornado dynamics \citep{GavrikovTaiurskii20,VaraksinRyzhkov23,RotunnoBluestein24}.

Tornadoes are fundamentally turbulent phenomena. To understand their structure and development, especially under realistic conditions, it is essential to view tornado dynamics in the context of turbulence. In \citet{PullinSaffman98}, turbulence is modeled using the three-dimensional incompressible Navier--Stokes equations, driven by initial conditions or external forcing, and vortex dynamics are recognized as central to turbulence evolution. \citet{MelanderHussain94} identified two key challenges in vortex dynamics: \textit{vortex core dynamics}, where the internal vorticity structure governs vortex evolution, and \textit{large-small scale interactions}, where small-scale vortices influence larger structures and vice versa. Related research has explored swirl effects in turbulence \citep{ShternHussain99}, sub-vortices within tornadoes \citep{FujitaEtAl70}, curvature-induced instabilities \citep{BlancoLe17}, vortex breakdown \citep{LiuEtAl18}, 
and horizontal vortex structures \citep{OliveiraEtAl19}.

Prior research by \citet{HsuNotsuYoneda16} on axisymmetric tornado-type flows in straight cylindrical domains has yielded valuable insights into vortex dynamics in idealized settings. Using hyperbolic inflow with and without swirl, they studied flow behavior near a saddle point. It is observed that only in the swirl case did the distance between the location where the maximum magnitude of velocity occurs and the $z$-axis change drastically at a specific time (called the turning point). Building on their approach using the Finite Element Method (FEM), we incorporate curvature into our simulations, considering that real tornadoes are often non-axisymmetric. While curvature significantly influences flow structures, energy distribution, and vortex behavior, its effects are not yet fully understood. By addressing this gap, our research enhances the understanding of vortex behavior in curved domains, with potential applications in engineering, fluid mechanics, and system design where curvature plays a critical role in flow dynamics.

In what follows, we give a summary of our result. Using FEM, we simulate incompressible Navier--Stokes equations with no-slip boundary and initial velocity with swirl on curved cylindrical domain. The swirl in the center of domain is imposed only in the initial time and with no external forcing applied thereafter. Observing simulations results, we aim to understand the effect of curvature. Observing velocity, we find that domain curvature exerts a greater influence than the initial velocity profile. It is shown that the location where the maximum magnitude of velocity occurs (hereafter, it will be abbreviated as maximum velocity location) for curved domain is gradually shifts outward, which higher curvature of the domain accelerates this outward motion. Investigation on low-pressure behavior reveals a small low-pressure region along $y$-axis that generated by flow dynamics. To discover more information, we also define the central curve. The low-pressure region and central curve also exhibit directional movement, in line with velocity observations. 

By calculating the kinetic energy and angular momentum from our simulations, we observe that the total kinetic energy and angular momentum are gradually decreasing. By dividing the domain into inner, middle, and outer regions based on distance to central curve, we observe energy transfer from inner to middle and then to outer region, with higher curvature accelerates this redistribution of momentum. To visualize vortex dynamics, we use Q-criterion \citep{JeongHussain95} with $Q\geq50$, $Q\geq250$ and  $Q\geq750$. The visualizations confirm the vortex outward movement, development of low-pressure region along $y$-axis, and energy transfer from inner to outer region. By defining primary and secondary vortices, we observe the emergence of a secondary vortex. The development of secondary vortex that appear in rear-left or southwest relative to the movement of primary vortex is in line with radar observation of real tornado (the El Reno, Oklahoma, Tornado on 31 May 2013) by \citet{BluesteinEtAl18}.

The remaining of the paper is organized as follows. Section \ref{sec:num_sim} presents the numerical simulation configuration, including the governing equations, numerical method, initial conditions, and domain settings. Section \ref{sec:results} presents the results of numerical simulations and discussion of velocity, pressure, kinetic energy, angular momentum, and vortex dynamics. Section \ref{sec:conclusion} provides the conclusions.

\section{Numerical Simulation}\label{sec:num_sim}

\subsection{Governing equations and numerical method}

Using FEM, we compute $(\boldsymbol{v},p) : \Omega \times (0,T) \rightarrow \mathbb{R}^3 \times \mathbb{R}$, where $ \boldsymbol{v} $ represents velocity and $ p $ represents pressure, for the incompressible Navier--Stokes equations with no-slip boundary condition:
\begin{equation}\label{eq:Navier--Stokes}
	\left\{ \begin{array}{lll}
		\partial_{t}\boldsymbol{v} + (\boldsymbol{v \cdot \nabla}) \boldsymbol{v} - \nu \Delta \boldsymbol{v} + \boldsymbol{\nabla} p = 0 & ,\boldsymbol{\nabla \cdot v} = 0 \ \mbox{ in } \Omega \times (0,T),\\[2pt]
		\boldsymbol{v} = \boldsymbol{0} \ \mbox{ on } \partial \Omega \times (0,T) & ,\boldsymbol{v} = \boldsymbol{v}^0 \ \mbox{ in } \Omega \mbox{ at } t=0.
	\end{array} \right.
\end{equation}
Here, $ \nu > 0 $ is the viscosity and $ \boldsymbol{v}^0 : \Omega \rightarrow \mathbb{R}^3 $ is the prescribed initial velocity.

The FEM used in this simulations is following prior research \citep{HsuNotsuYoneda16} which is the stabilized Lagrange--Galerkin scheme, cf.~\citet{NotsuTabata08,NotsuTabata15a,NotsuTabata15b}, to find a pair of piecewise linear functions $(\boldsymbol{v}_h^{k},p_h^{k})$, approximation of $(\boldsymbol{v},p)$ at $t = k\tau$, in a {\it strong}-representation:
\begin{align*}{}
    \dfrac{1}{\tau} \Bigl[ \boldsymbol{v}_h^{k} (\boldsymbol{x}) - \boldsymbol{v}_h^{k-1} \Bigl( \boldsymbol{x} - \boldsymbol{v}_h^{k-1} (\boldsymbol{x}) \tau \Bigr) \Bigr] - \nu \Delta\boldsymbol{v}_h^{k} + \boldsymbol{\nabla} p_h^k & = 0, \\[2pt]
    \boldsymbol{\nabla \cdot v}_h^{k} - \delta^{\rm s}_0\, h^2 \Delta p_h^k & = 0,
\end{align*}
for $k=1,2,\dots, N_T$, where $\tau, h, \delta^{\rm s}_0 > 0$, and $N_T \coloneq \lfloor T/\tau \rfloor \in\mathbb{N}$ are a time-step size, a mesh size, a stabilization parameter, and a total number of time steps, respectively.
The first term in the first equation of the scheme is based on the idea of the method of characteristics, i.e.,
\[
    \Bigl[ \partial_{t}\boldsymbol{v} + (\boldsymbol{v \cdot \nabla}) \boldsymbol{v} \Bigr](\boldsymbol{x}, k\tau) \approx
    \dfrac{1}{\tau} \Bigl[ \boldsymbol{v}^{k} (\boldsymbol{x}) - \boldsymbol{v}^{k-1} \Bigl( \boldsymbol{x} - \boldsymbol{v}^{k-1} (\boldsymbol{x}) \tau \Bigr) \Bigr],
\]
with the notation $ \boldsymbol{v}^{k} (\boldsymbol{x}) \coloneq \boldsymbol{v}(\boldsymbol{x}, k\tau)$ for the approximation of the material derivative,
and the second term $-\delta^{\rm s}_0\, h^2 \Delta p_h^k$ in the second equation is the pressure-stabilization introduced by \citet{BrePit-1984}. It is worth noting that the scheme converges to the exact solution with accuracy of order $O(\tau + h)$ in $L^\infty(0,T; H^1(\Omega;\mathbb{R}^3))$ and of order $O(\tau + h^2)$ in $L^\infty(0,T; L^2(\Omega;\mathbb{R}^3))$ under some conditions, see~\citet{NotsuTabata15b}. In this paper, we chose $\tau = 1.25 \times 10^{-2}$ and $\delta^{\rm s}_0=1$.

\subsection{Domain and initial configurations}

For simulations, we consider swirl initial velocity 
\begin{equation} \label{eq:initial}
    \boldsymbol{v}^0 \coloneq u_r \boldsymbol{e}_r + u_\theta \boldsymbol{e}_\theta + u_z \boldsymbol{e}_z,
\end{equation}
where $ \boldsymbol{e}_r = (1/ \sqrt{x^2+y^2})(x,y,0)$, $ \boldsymbol{e}_\theta = (1/ \sqrt{x^2+y^2})(-y,x,0) $, $ \boldsymbol{e}_z = (0,0,1) $
with:
\begin{align*}
    u_{z} &= \psi(r,\epsilon_{1},-\beta_{1})\psi(z,\epsilon_{2},-\beta_{2}),\\
    \rho &= \psi(r,\epsilon_{3},-\beta_{3})\psi(z,\epsilon_{4},\beta_{4}),\\
    u_{\theta} &= \psi(r,\epsilon_{5},-\beta_{5})\psi(z,\epsilon_{6},-\beta_{6}),\\
    u_{r} &= \operatorname{sign}(z)\rho u_{z},
\end{align*}
for $ \psi(a,\epsilon,\sigma) = (a^{2}+\epsilon)^{\sigma} $.
The constants $ \epsilon_{i} \mbox{ and } \beta_{i} \ (i=1,\dots,6) $ are set to $1$. By this setting, the initial magnitude of velocity $|\boldsymbol{v}|$ at $(x,y,z)=(0,0,0)$ (called center of the initial velocity, which close to lower boundary but not on boundary) was larger than at other places (see figure~\ref{fig:initial_velocity_swirl}). Note that the configuration for initial velocity and straight cylindrical domain are similar to those used in \citet{HsuNotsuYoneda16} where it numerically satisfy the divergence-free and no-slip boundary condition after the first time step, not at the initial step.

\begin{figure}[!h]
\centerline{\includegraphics[width=1\textwidth]{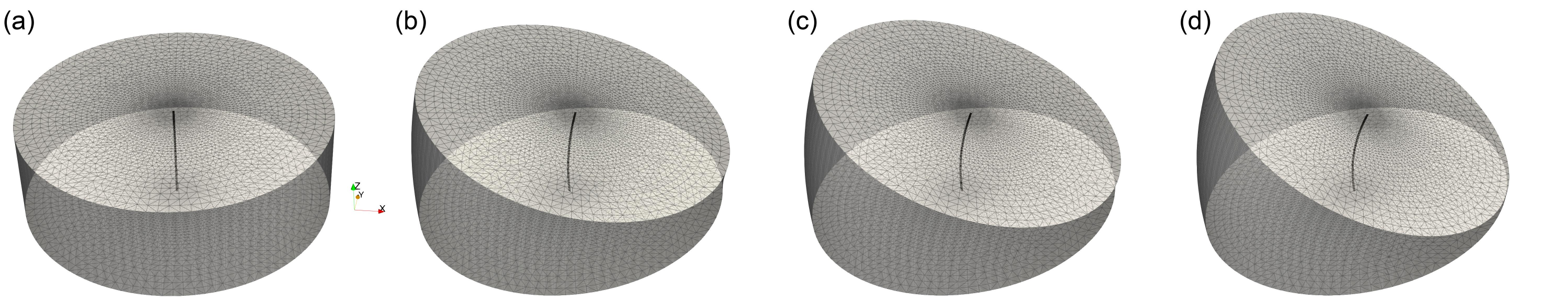}}
	\caption{3D visualization of the domain shapes of (a) the straight cylindrical domain $ \Omega_{S} $ and the curved cylindrical domain $ \Omega_{C} $  with its geometric axis (black curve in center): (b) $ \Omega_{C} $ with $ R=2.0 $, (c) $ \Omega_{C} $ with $ R=1.5 $, and (d) $ \Omega_{C} $ with $ R=1.1 $.}
\label{fig:3d_domain}
\end{figure}

\begin{figure}[!h]
\centerline{\includegraphics[width=0.8\textwidth]{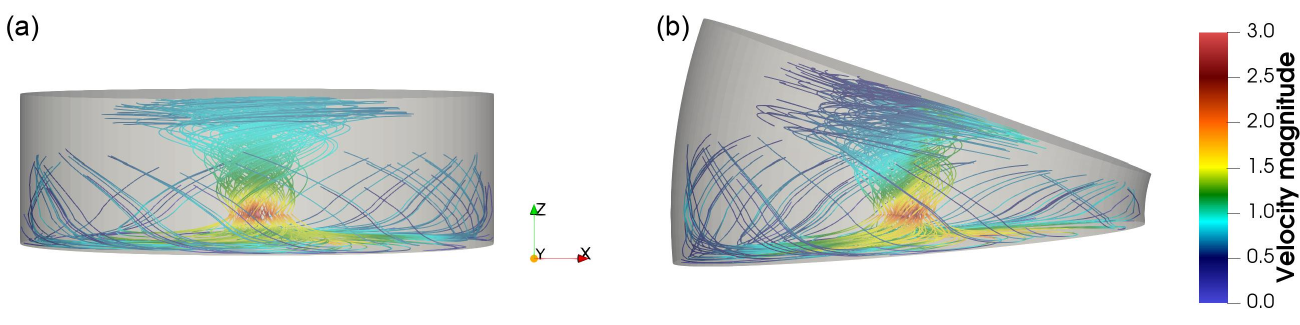}}
	\caption{3D swirl initial velocity streamlines in (a) the straight cylindrical domain $ \Omega_{S} $ and (b) the curved cylindrical domain $ \Omega_{C} $ with $ R=1.5 $.}
\label{fig:initial_velocity_swirl}
\end{figure}

The domains considered in this study are straight cylindrical domain and curved cylindrical domains. The straight cylindrical domain (figure~\ref{fig:3d_domain}a), $ \Omega_{S} $, is defined as  
\begin{equation}\label{eq:straight_domain}
	\Omega_{S} \coloneq \{\boldsymbol{x}_{S} = (x,y,z) \in \mathbb{R}^{3} \mid -a \leq z \leq 4a ,  \ r(\boldsymbol{x}_{S}) \leq r_{\max}\}
\end{equation}
where $ r_{\max} \in \mathbb{R} $  is the maximum radius, and $ r,a > 0 $. To model the curved cylindrical domain (figure~\ref{fig:3d_domain}b-d) $ ,\Omega_{C} $ (conceptualized as a toroidal segment), we apply the mapping $ T_{S,R} : \Omega_{S} \rightarrow \Omega_{C} $, given by
\begin{equation}\label{eq:transform_domain}
	T_{S,R}(\boldsymbol{x}_{S}) \coloneq \Big(R-(R-x)\cos(\dfrac{z}{R}) ,  y ,  (R-x)\sin(\dfrac{z}{R})\Big)
\end{equation}
where $ \boldsymbol{x}_{S} \in \Omega_{S} $ and $ R>r_{\max} $ is torus radius. The curved cylindrical domain is defined as
\begin{equation}\label{eq:curved_domain}
	\Omega_{C} \coloneq \{ \boldsymbol{x}_{C} \in \mathbb{R}^{3} \mid \boldsymbol{x}_{C} = T_{S,R}(\boldsymbol{x}_{S}) , \ r_{\max}\leq R\}.
\end{equation}

For all simulations, we set $r_{\max}=1$, $a=0.125$, and Reynolds number $ \Rey=10^4 $. To ensure turbulent flow while maintaining reasonable computational efficiency, only a single representative Reynolds number is used. This approach allows us to concentrate on the influence of curvature variation, which is the primary focus of this study. We consider three torus radii $R=$ 2, 1.5, and 1.1, where smaller values of $ R $ correspond to larger curvature. As our curved cylindrical domain are conceptualized as a toroidal segment, by definition of curvature in toroidal $\delta = r_{\max} / R$, then our simulations consider curvature $\delta=$ 0.5, 0.667, and 0.909. In this study, we distinguish between the geometric axis (figure~\ref{fig:3d_domain}), which refers to the fixed and time-independent at the center of the domain, and the central curve (figure~\ref{fig:3d_low_pressure_region}), which is flow-dependent and may evolve over time.

\begin{figure}	[!h]\centerline{\includegraphics[width=1\textwidth]{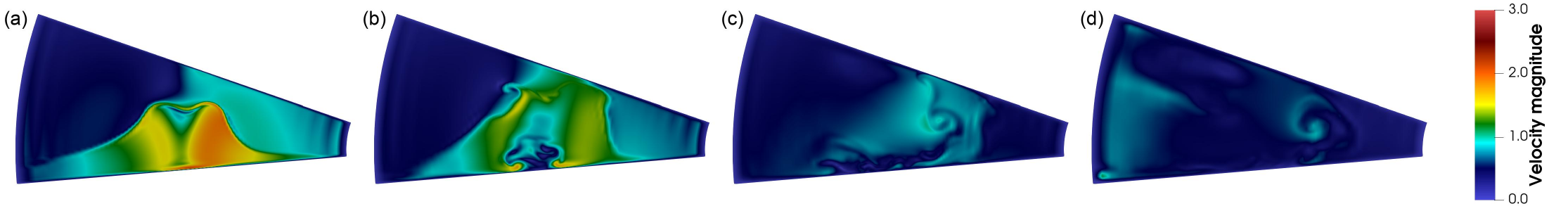}}
	\caption{Time evolution of the $xz$-plane slice of velocity magnitude $ |\boldsymbol{v}| $ for the curved cylindrical domain $ \Omega_{C} $ ($ R=1.5 $) at (a) $ t = 0.3 $, (b) $ t = 0.8 $, (c) $ t = 1.9 $, and (d) $ t = 2.9 $.}
\label{fig:2d_slice_velocity}
\end{figure}

\section{Results and Discussion}\label{sec:results}

The simulations results are presented to understand the effect of curvature of the domain by observing the velocity, pressure, kinetic energy, angular momentum, and vortex dynamics.

\subsection{Velocity observations}

\begin{figure}[!h]
	\centerline{\includegraphics[width=1\textwidth]{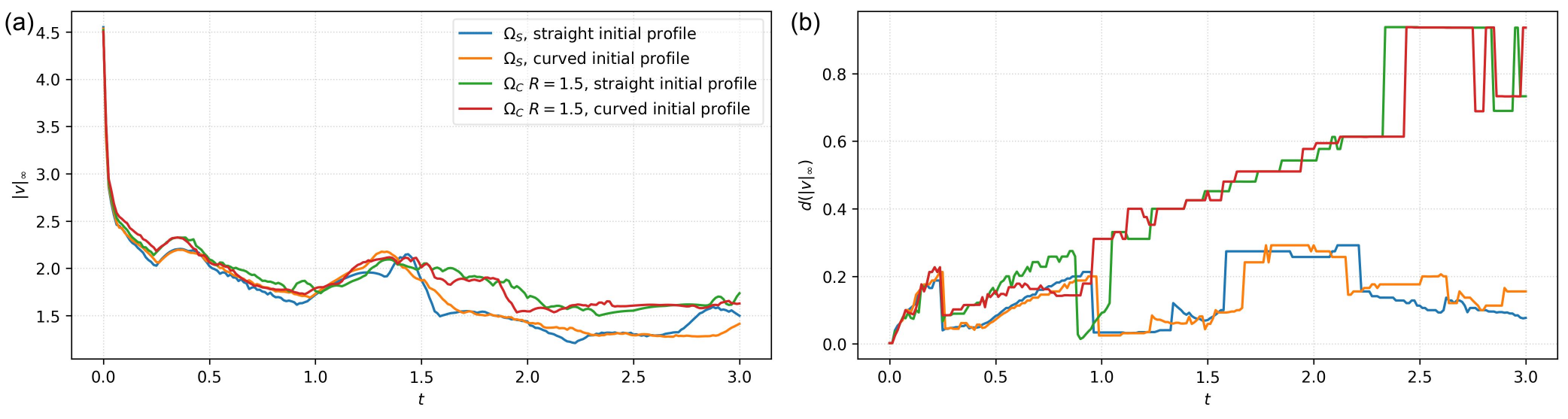}}
	\caption{Time evolution of (a) the maximum magnitude of velocity $ |\boldsymbol{v}|_{\infty} $ and (b) its distance to the geometric axis $ d(|\boldsymbol{v}|_{\infty}) $, for combinations of domain (straight cylindrical domain $ \Omega_{S} $ and curved cylindrical domain $ \Omega_{C} $ with $ R=1.5 $) and initial velocity profile (straight and curved).}
\label{fig:graph_velocity_of_combination}
\end{figure}

For curved cylindrical domain $ \Omega_{C}$ $(R=1.5)$, the initial velocity streamlines are shown in figure~\ref{fig:initial_velocity_swirl}b and the evolution of velocity magnitude in $xz$-plane shown in figure~\ref{fig:2d_slice_velocity}. By default, the swirl initial velocity profile for the simulation is designed to align the domain shape (figure~\ref{fig:initial_velocity_swirl}). However, to investigate the influence of the initial velocity profile on the simulation, we compare results by pairing domain shapes ($ \Omega_{S} $ and $ \Omega_{C} $) with initial velocity profiles (straight and curved). We analyze the maximum magnitude of velocity, denoted by $ |\boldsymbol{v}|_{\infty}(t) =|\boldsymbol{v}|_{\infty} $, and its distance from the geometric axis or maximum velocity location, $ d(|\boldsymbol{v}|_{\infty}) $. The results, shown in figure~\ref{fig:graph_velocity_of_combination}, show that the simulation outcomes (both $ |\boldsymbol{v}|_{\infty} $ and $ d(|\boldsymbol{v}|_{\infty}) $) are similar for the same domain shape, regardless of the initial velocity profile. This suggests that the domain geometry has a greater influence on the simulation results than the choice of initial velocity profile. 

\begin{figure}[!h]
    \centerline{\includegraphics[width=1\textwidth]{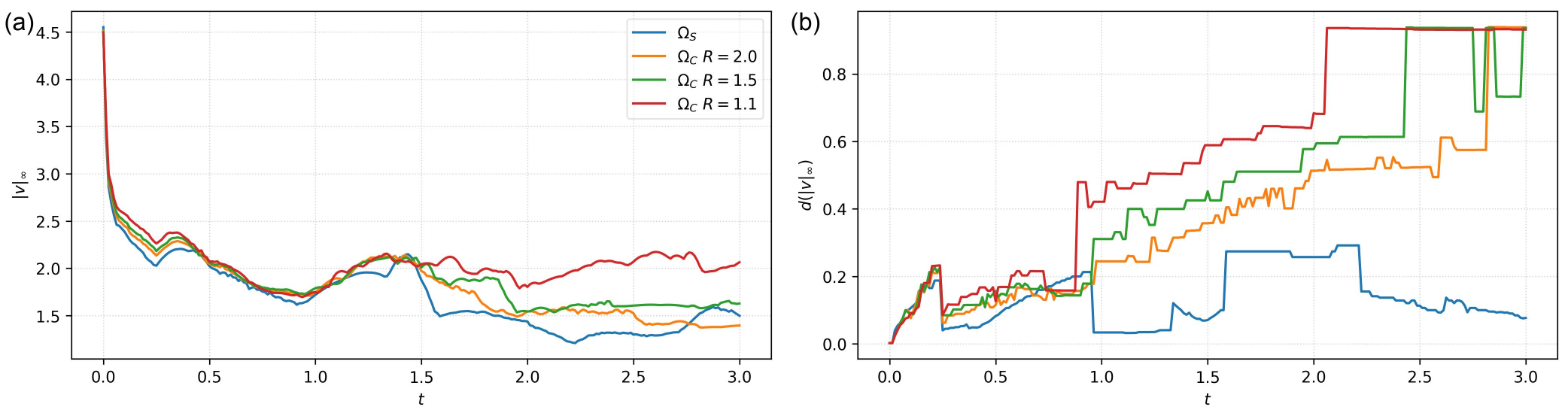}}
	\caption{Time evolution of (a) the maximum magnitude of velocity $ |\boldsymbol{v}|_{\infty} $ and (b) its distance to the geometric axis $ d(|\boldsymbol{v}|_{\infty}) $, for the straight cylindrical domain $ \Omega_{S} $ and curved cylindrical domain $ \Omega_{C} $ with $ R=2 $, $R=1.5$, and $R=1.1$.}
\label{fig:graph_velocity_of_R}
\end{figure}

The influence of domain shape is investigated through simulations of the curved cylindrical domain $ \Omega_{C} $ with varying curvature. Results in figure~\ref{fig:graph_velocity_of_R} show that in $ \Omega_{C} $, maximum velocity location moves progressively farther from the geometric axis, unlike in $ \Omega_{S} $, where it remains close. This suggests that curvature of the domain causes the maximum velocity location to move outward. Furthermore, higher curvature (smaller $ R $) results in larger $ d(|\boldsymbol{v}|_{\infty}) $, indicating that the maximum velocity location shifts farther from the geometric axis.

\begin{figure}[!h]
	\centerline{\includegraphics[width=1\textwidth]{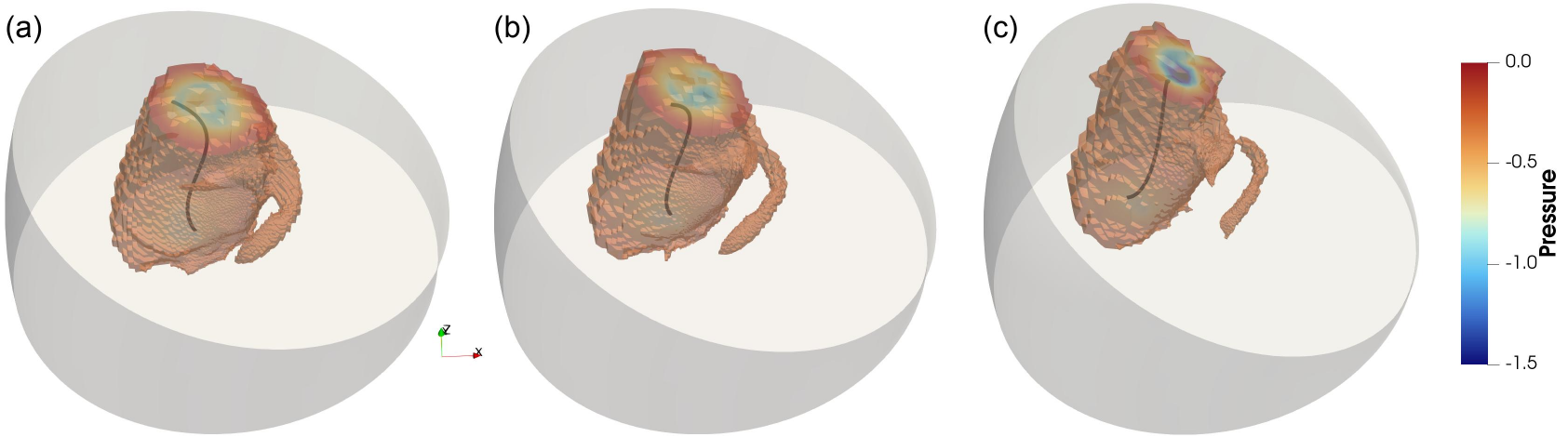}}
	\caption{3D visualization of the low-pressure region with the central curve (black) at $t=1.9$ in curved cylindrical domains $ \Omega_{C} $ for torus radii (a) $ R=2 $, (b) $ R=1.5 $, and (c) $ R=1.1 $.}
\label{fig:3d_low_pressure_region}
\end{figure}

\subsection{Pressure observations}

To further understand the flow behavior in curved cylindrical domains $ \Omega_{C} $, we analyze low-pressure regions alongside $ |\boldsymbol{v}|_{\infty} $ and $ d(|\boldsymbol{v}|_{\infty}) $. Initially, the low-pressure region forms a large, tube-like structure aligned with the $ z $-axis, consistent with the initial conditions. By $ t=1.9 $, a smaller low-pressure region develops along the $ y $-axis (figure~\ref{fig:3d_low_pressure_region}), indicating that this phenomenon arises from the flow dynamics rather than initial conditions. In figure~\ref{fig:3d_low_pressure_region}, it can also be observed that the low-pressure region is not centered within the domain. This indicates that the low-pressure region shift further from geometric axis, aligning with the displacement of the maximum velocity location.

Consider a set of pressure points $ M_{t} \coloneq \{\boldsymbol{P}_{j}\}_{j
=1}^{N_{p}} $, where $ N_p $ is the number of nodes, $ t $ is the time index, and $ \boldsymbol{P} $ represents the location vector of the pressure point. Using the definitions in \eqref{eq:straight_domain}, \eqref{eq:transform_domain}, and \eqref{eq:curved_domain}, the cross section plane-$l$ in straight cylindrical domain is defined as $ \Pi ^{S}_{l} = \{ (x,y,z)\in \Omega_S \mid z = z_l , \ x^2+y^2 \leq r_{\max} \} $ and for the curved domain $ \Pi ^{C}_{l} = T_{S,R}(\Pi ^S_l) $ for $ l = 0, \dots, N_{l} $ (where $ N_{l}=100 $). Note that $z_l \coloneq z_{\min} + l \Delta z \in [z_{\min}, z_{\max}]$, where $z_{\min} \coloneq -a$, $\Delta z \coloneq 5a/N_l$ and $z_{\max} \coloneq 4a$. The minimum pressure point in each plane-$l$ is given by:
\[
    \boldsymbol{P}_{\min,l} \coloneq \operatorname*{argmin}_{\boldsymbol{P} \in M_{t,l}} \{ p_{t,l}(\boldsymbol{P})\},
\]
with $p_{t,l}$ is the pressure at time~$t$ on plane-$l$, where $ p_{t}(\boldsymbol{P}_{\min,l}) \leq p_{t}(\boldsymbol{P}_{m}) $ for all $ m
$ in $ M_{t,l} \coloneq M_{t} \cap \Pi_{l} $ for $\Pi_{l} = \Pi_{l}^{S}$ or $\Pi_{l}^C$. We introduce the central curve $ \{ \boldsymbol{C}(\xi) \in \mathbb{R}^3 \mid \xi \in [z_{\min},z_{\max}] \} $. For a cubic polynomial representation of the central curve: 
\[
    \boldsymbol{C}(\xi) = \begin{bmatrix} c_{0}^{(1)} + c_{1}^{(1)} \xi + c_{2}^{(1)} \xi^{2} + c_{3}^{(1)} \xi^{3} \\ c_{0}^{(2)} + c_{1}^{(2)} \xi + c_{2}^{(2)} \xi^{2} + c_{3}^{(2)} \xi^{3}\\ c_{0}^{(3)} + c_{1}^{(3)} \xi + c_{2}^{(3)} \xi^{2} + c_{3}^{(3)} \xi^{3} \end{bmatrix}, 
\] 
we minimize the function $ J( \boldsymbol{c} = [c_{0}^{(1)}, \ldots, c_{3}^{(3)}] ) \coloneq \thalf (\sum_{l} |\boldsymbol{C} (\xi_{l}) - \boldsymbol{P}_{\min,l}|^{2}) $ to determine the coefficients $ \boldsymbol{c} $. Substituting $ \boldsymbol{c} $ and the parameter $ \xi $ into $\boldsymbol{C}(\xi)$ yields the central curve. 

\begin{figure}[!h]
	\centerline{\includegraphics[width=1\textwidth]{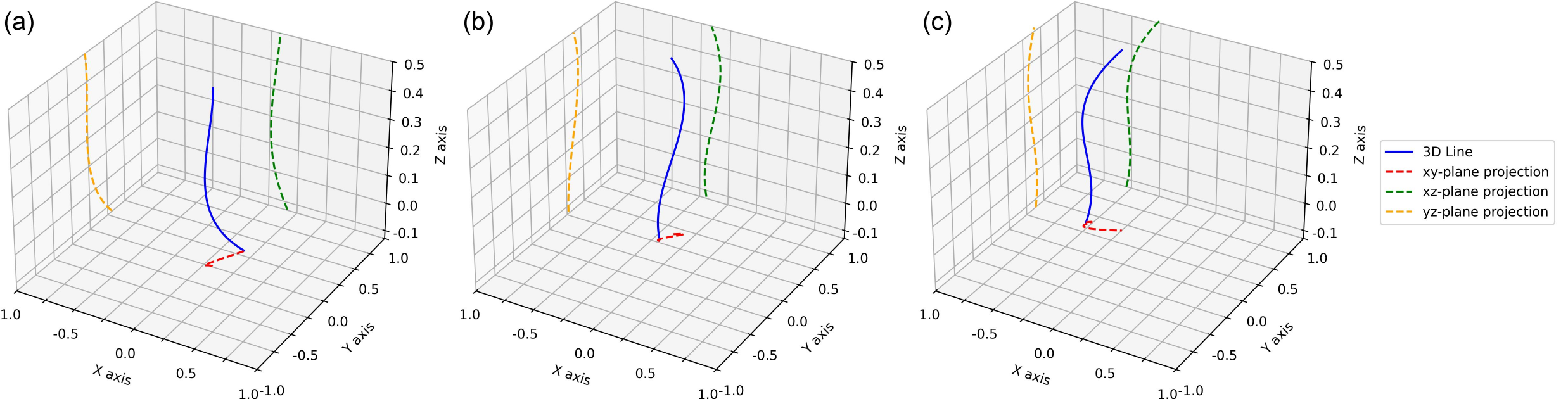}}
	\caption{Time evolution of the 3D central curve (blue) with its projection on 2D planes for curved cylindrical domain $\Omega_{C}$ with $R=1.5$ at (a) $t=0$, (b) $t=1.9$, and (c) $t=3$.}
\label{fig:3d_central_line_projection}
\end{figure}

Detailed projections of central curve for the curved cylindrical domain $\Omega_{C}$ with $R=1.5$ are shown in figure~\ref{fig:3d_central_line_projection}, revealing a shift toward the negative $ x $-axis and positive $ y $-axis. The average movement speeds in the same directions for $ R = 2$ , $1.5$, and $1.1$, are $0.366$, $0.418$, and $0.49$, respectively, indicating that smaller $ R $ results in faster outward motion. This trend aligns with the earlier observation of $ d(|\boldsymbol{v}|_{\infty}) $ (figure~\ref{fig:graph_velocity_of_R}b) and further supports the connection between low-pressure regions and high-velocity regions, reinforcing their role in defining vortex structures.

\subsection{Energy transfer via kinetic energy and angular momentum}

\begin{figure}[!h]
	\centerline{\includegraphics[width=1\textwidth]{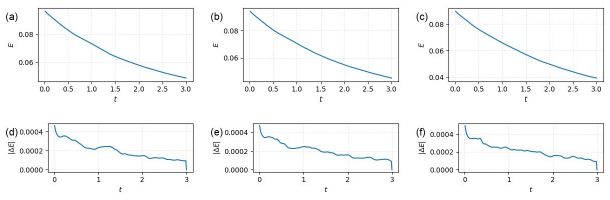}}
	\caption{Time evolution of total energy kinetic $E$ and its absolute change $|\Delta E|$, for curved cylindrical domain $\Omega_{C}$ for torus radii $ R=2 $ (a,d), $ R=1.5 $ (b,e), and $ R=1.1 $ (c,f).}
\label{fig:graph_total_energy}
\end{figure}

To understand energy transfer mechanisms in the flow, we analyze both kinetic energy and angular momentum. Kinetic energy ($ E $) quantifies the total fluid motion, where higher values indicate more energetic flow structures:
\begin{equation}
    E = \int_{\Omega_{C}} \frac{1}{2} |\boldsymbol{v}|^{2} \ d\Omega_{C} \approx \frac{1}{2} 
    \sum^{N_e} |\boldsymbol{v}_{e}|^{2} \ V_{e}
\end{equation}
with total number of elements $N_e$ where $ \boldsymbol{v}_{e} $ is the velocity vector and $ V_{e} $ the volume of the tetrahedral element. Figure~\ref{fig:graph_total_energy} shows a gradual energy decay in all simulations with various torus radii. The gradual loss of total kinetic energy suggests that the turbulence is decaying over time.

\begin{figure}[!h]
	\centerline{\includegraphics[width=1\textwidth]{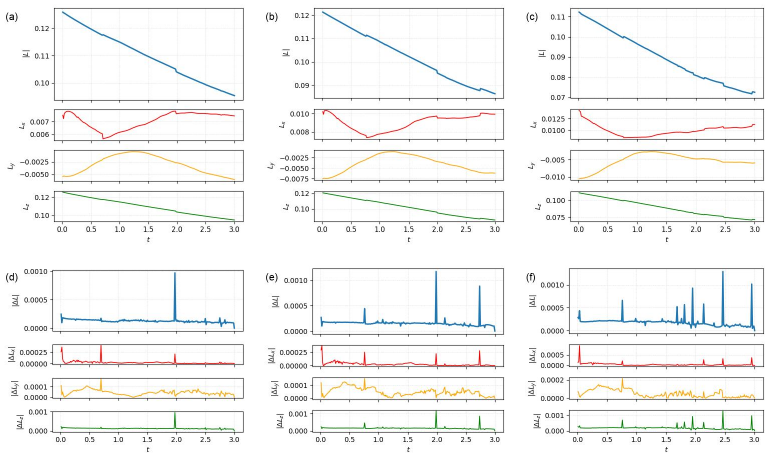}}
	\caption{Time evolution of total angular momentum $ \boldsymbol{L} $ and its absolute change $|\Delta \boldsymbol{L}|$ per time step, for curved cylindrical domain $\Omega_{C}$ for torus radii $ R=2 $ (a,d), $ R=1.5 $ (b,e), and $ R=1.1 $ (c,f) with value of the total (blue), x-component (red), y-component (yellow), and z-component (green).}
\label{fig:graph_total_angular_momentum}
\end{figure}

Angular momentum ($ \boldsymbol{L} $) characterizes the rotational motion:
\begin{equation}
	\boldsymbol{L} = \int_{\Omega_{C}} (\boldsymbol{r}_{e} \times \boldsymbol{v}_{e}) \ d\Omega_{C} \approx \sum^{N_e} (\boldsymbol{r}_{e} \times \boldsymbol{v}_{e}) \ V_{e}
\end{equation}
where $ \boldsymbol{r}_{e} = \boldsymbol{x}_{e} - \boldsymbol{C}_{e} $ is the relative position vector of element $e$, defined as the displacement from the nearest point on central curve $\boldsymbol{C}_e$ to the element centeroid $\boldsymbol{x}_e$. While angular momentum would be conserved in an ideal system, the no-slip boundary conditions cause $|\boldsymbol{L}|$ to decrease over time (figure~\ref{fig:graph_total_angular_momentum}a-c).

Figure~\ref{fig:graph_total_angular_momentum} shows that the total angular momentum is primarily dominated by its $z$-component, $ L_z $ (green), which is expected due to our initial condition that induce tornado-type of flow. It is also observed that the total angular momentum gradually decays like the kinetic energy. However, this dissipation is not directly coupled to the kinetic energy dissipation rate. Observed spikes in $|\Delta \boldsymbol{L}|$ (figure~\ref{fig:graph_total_angular_momentum}d-f) correspond to brief momentum transfer events or reorganization, which occur more frequently at higher curvature $ (R=1.1) $, suggesting stronger dynamic effects. Notably, distinct peaks appear around $t=0.8$ and $t=1.9$ across all simulation with various torus radius. The sudden reorganization of rotation results in peaks in $|\Delta \boldsymbol{L}|$, but the $|\Delta E|$ (which includes both rotational and non-rotational motion) may not exhibit spikes due to internal redistribution rather than total dissipation. These angular momentum peaks, which are not observed in the kinetic energy, may reflect localized vortex interactions.

To investigate flow dynamics, we divide the domain into three regions based on the radial distance $ |r_{e}| $ from central curve: the inner region $  r_{e,1} $ for $ |r_{e}| \leq 0.15 $, the middle region $ r_{e,2} $ for $ 0.15 < |r_{e}| \leq 0.4 $, and the outer region $ r_{e,3} $ for $ 0.4 < |r_{e}| \leq 0.7 $. This division, guided by low-pressure regions and vortex identification (via Q-criterion visualization), enables clearer analysis of energy distribution and transfer between regions. Notably, the inner region has the smallest volume, the middle region is larger, and the outer region has the largest volume.

\begin{figure}[!h]
	\centerline{\includegraphics[width=1\textwidth]{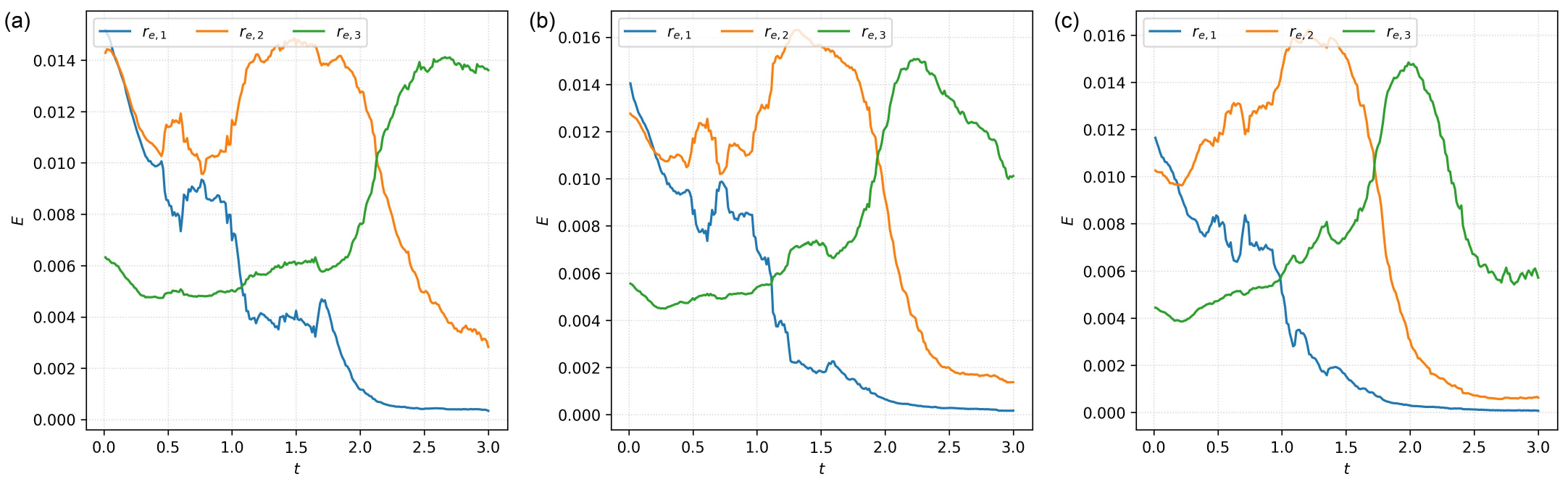}}
	\caption{Time evolution of the energy kinetic $E$ for the inner region $  r_{e,1} $ (blue), middle region $ r_{e,2} $ (orange), and outer region $ r_{e,3} $ (green) in the curved cylindrical domain $\Omega_{C}$ with (a) $ R=2 $, (b) $ R=1.5 $, and (c) $ R=1.1 $.}
\label{fig:graph_energy_exchange}
\end{figure}

\begin{figure}[!h]
	\centerline{\includegraphics[width=1\textwidth]{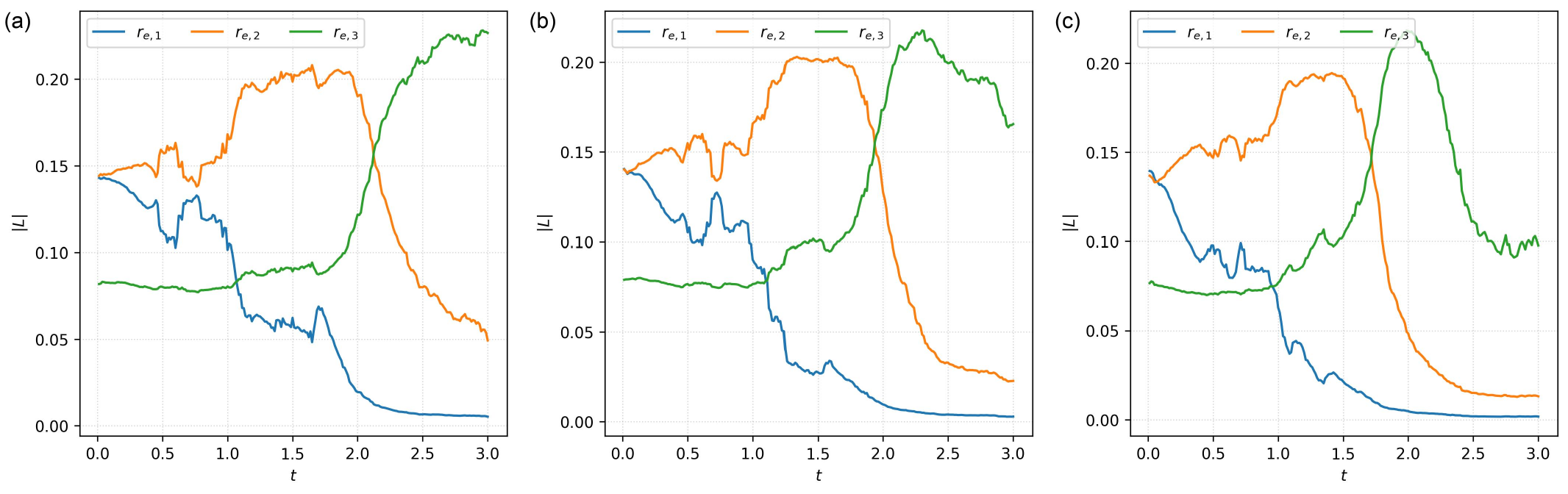}}
	\caption{Time evolution of the total angular momentum $|\boldsymbol{L}|$ for the inner region $  r_{e,1} $ (blue), middle region $ r_{e,2} $ (orange), and outer region $ r_{e,3} $ (green) in the curved cylindrical domain $\Omega_{C}$ with (a) $ R=2 $, (b) $ R=1.5 $, and (c) $ R=1.1 $.}
\label{fig:graph_angular_momentum_exchange}
\end{figure}

The energy and angular momentum dynamics are revealed in figures \ref{fig:graph_energy_exchange} and \ref{fig:graph_angular_momentum_exchange}, which show distinct fluctuation patterns, regional concentrations, and correlations with vortex activity. Given that only the initial velocity is prescribed (no external force) and angular momentum is conserved within the constrained domain, energy decreases in one region necessarily correspond to increases elsewhere. This conservation principle illuminates the rotational intensity and distribution, demonstrating how vortex structures facilitate energy redistribution.

In the inner region ($r_{e,1}$), both kinetic energy ($E$) and angular momentum ($|\boldsymbol{L}|$) exhibit a rapid decline for higher curvature ($R=1.1$), while decreasing more gradually for lower curvature. The middle region ($r_{e,2}$) shows a distinctive sharp peak followed by decay, particularly pronounced at higher curvature, indicating its role as an intermediate zone for energy and momentum transfer between the inner and outer regions. Meanwhile, the outer region ($r_{e,3}$) shows steady increases in both $E$ and $|\boldsymbol{L}|$, with their parallel evolution reflecting coherent vortex-mediated transport. Notably, the high curvature case ($R=1.1$) exhibits a sharp energy drop in the outer region around $t=1.9$, whereas for $R=1.5$ this occurs slightly later at around $t=2.3$. In contrast, the lowest curvature case ($R=2$) maintains stable outer-region momentum throughout the simulation period. These sudden drops in the outer region may be influenced by the flow interacting with the no-slip boundary wall.

The transfer of angular momentum across regions is shown in figure~\ref{fig:graph_angular_momentum_exchange}. Transfer from the inner region ($r_{e,1}$) to the middle region ($r_{e,2}$) happens around $t=0.8$ for all torus radii. This is indicated by the blue angular momentum profile decreasing while the orange angular momentum profile increases. Then, the transform from the middle region ($r_{e,2}$) to the outer region ($r_{e,3}$) occurs around $t=1.9$ for $R=1.5$, earlier for $R=1.1$ (higher curvature), and later for $R=2$ (lower curvature). This is indicated by the orange angular momentum profile decreasing while the green angular momentum profile increases.

The influence of domain's curvature manifests clearly in these dynamics. Higher curvature ($R=1.1$) accelerates both the central curve movement and the momentum redistribution from inner to outer regions. Conversely, lower curvature ($R=2$) promotes flow stability, resulting in weaker momentum exchange and more gradual evolution of the central curve. These observations collectively demonstrate how increased curvature enhances transport efficiency while decreasing curvature leads to more uniform, less dynamic flow patterns.

\subsection{Vortex identification using the Q-criterion}

Defining vortex structures in three-dimensional flows remains a critical issue, as a universally accepted definition of a vortex is still lacking. Loosely speaking, a vortex refers to a region where vorticity levels are higher than those in the surrounding flow. Therefore, in this study, we distinguish between primary and secondary vortices based on their origin and observed behavior. The primary vortex refers to the dominant swirling structure initialized in our simulations, which maintains its coherence as the largest and most intense flow feature. Secondary vortices are non-initialized structures that emerge around the primary vortex due to flow interactions, typically appearing as smaller, shorter-lived features, oriented perpendicular to the primary vortex and located at mid-level heights.

To objectively identify these features, we employ the Q-criterion method \citep{JeongHussain95}, which evaluates the local balance between rotational and straining motions through the velocity gradient tensor. It is defined as:
\begin{equation}
	Q = \dfrac{1}{2} \Big( |\mathbf{W}|^{2} - |\mathbf{S}|^{2} \Big) = \dfrac{1}{2} \Big( \sum_{ij}\mathbf{W}_{ij}^{2} - \sum_{ij} \mathbf{S}_{ij}^{2} \Big)
\end{equation}
where $ \mathbf{W} = \thalf \big( \nabla \boldsymbol{v} - (\nabla \boldsymbol{v})^{T}\big) $ is the antisymmetric part (representing rotational motion) and $ \mathbf{S} = \thalf \big( \nabla \boldsymbol{v} + (\nabla \boldsymbol{v})^{T}\big) $ is the symmetric part (representing strain rate). A vortex is identified when rotational motion dominates over strain, indicated by $ Q > 0 $.

\begin{figure}[!h]
	\centerline{\includegraphics[width=0.9\textwidth]{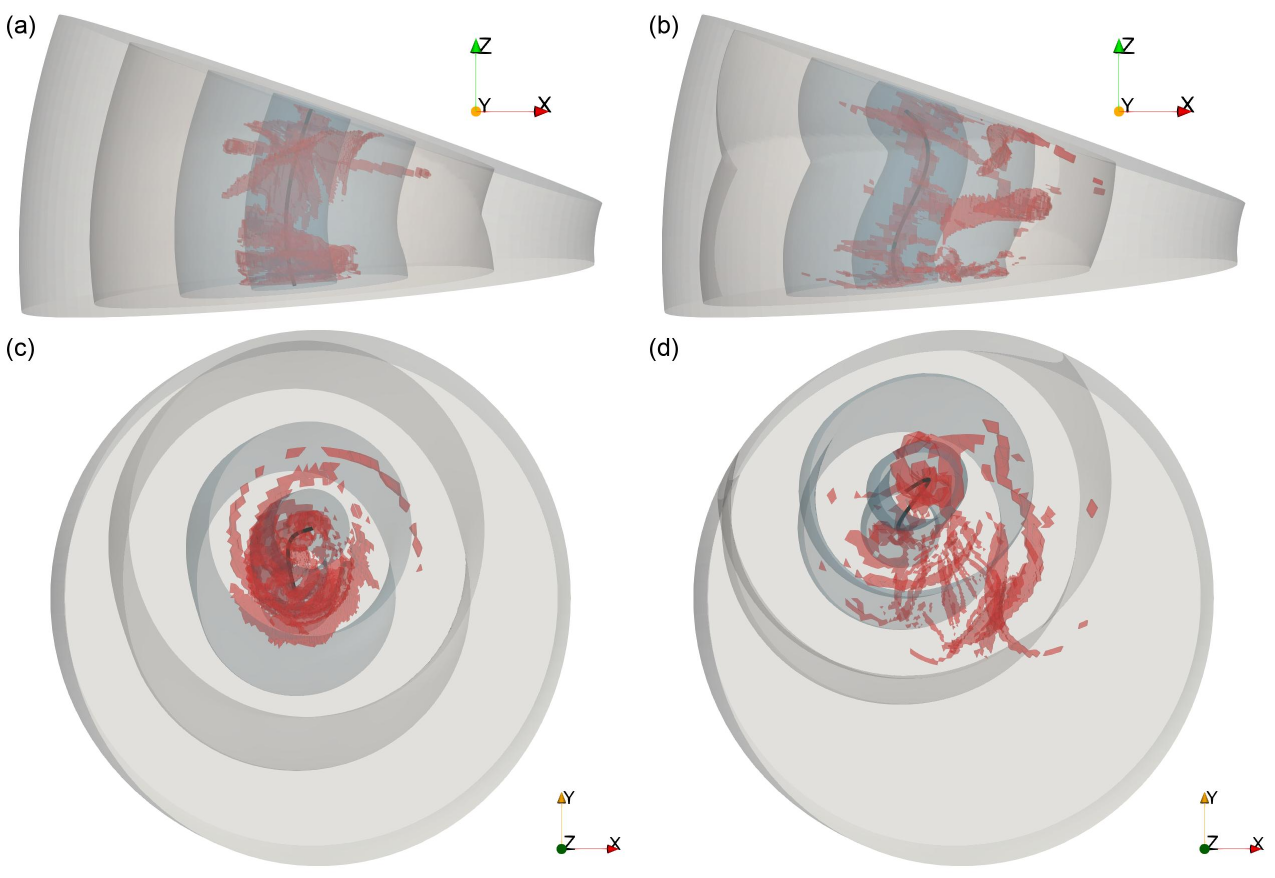}}
	\caption{3D vortex in the curved cylindrical domain $ \Omega_{C} $ ($ R=1.5 $) using Q-criterion for $Q \geq 50$ (red) at $t=0.8$ (a,c) and $t=1.9$ (b,d) with front-(a,b) and top-views(c,d). central curve (black), inner $  r_{e,1} $ (blue), middle $ r_{e,2} $ (light blue), and outer $ r_{e,3} $ (gray) regions are shown.}
\label{fig:3d_vortex_Qcriterion_50}
\end{figure}

Visualizing the Q-criterion for $ Q \geq 50$ in the curved cylindrical domain $ \Omega_{C} $ with $ R=1.5 $, we obtain figure~\ref{fig:3d_vortex_Qcriterion_50}, shown for $ t=0.8 $ and $ t=1.9 $. Comparing these results with the angular momentum exchange trends in figure~\ref{fig:graph_angular_momentum_exchange}b, we observe a consistent pattern of vortex transport. At $ t=0.8 $, vortices move outward, reducing the angular momentum in the inner region and increasing it in the middle region. By $ t=1.9 $, angular momentum shifts further outward, accumulating in the outer region, highlighting rotational energy redistribution. Vortex structures align with low-pressure regions (figure~\ref{fig:3d_low_pressure_region}b), particularly the vortex along the $ y $-axis, also captured in the velocity magnitude slice (figure~\ref{fig:2d_slice_velocity}c). As studied by \citet{GotoEtAl17}, the formation of this secondary vortex 
ensures angular momentum conservation and energy transfer to smaller scales, which associated with vortex creation events, consistent with localized energy cascade theory.

\begin{figure}[!h]
	\centerline{\includegraphics[width=1\textwidth]{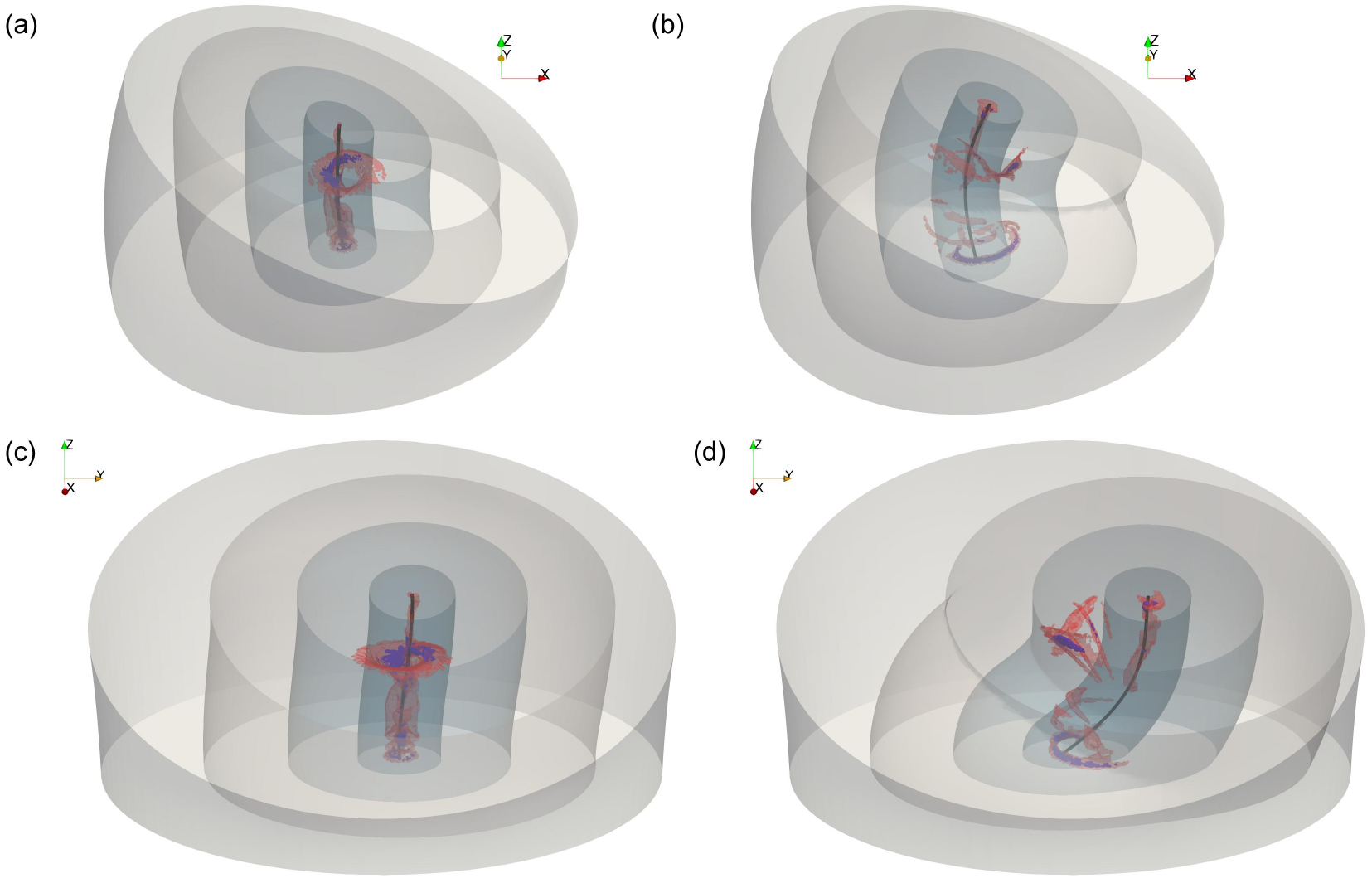}}
	\caption{3D vortex in the curved cylindrical domain $ \Omega_{C} $ ($ R=1.5 $) using Q-criterion for $Q \geq 250$ (red) and $Q \geq 750$ (dark blue) at $t=0.3$ (a,c) and $t=0.9$ (b,d) with front-(a,b) and side-views(c,d). central curve (black), inner $  r_{e,1} $ (blue), middle $ r_{e,2} $ (light blue), and outer $ r_{e,3} $ (gray) regions are shown.}
\label{fig:3d_vortex_Qcriterion_250_750}
\end{figure}

The Q-criterion analysis clearly distinguishes vortex structures by intensity, with strong vortices (blue, $ Q\geq750 $) and weaker vortices (red, $ Q\geq250 $) visible in figure~\ref{fig:3d_vortex_Qcriterion_250_750}. At $t=0.3$, the primary vortex aligns with the central curve at the pressure minimum. By $t=0.9$, two significant developments occur: (i) a secondary vortex forms perpendicular to the primary vortex, and (ii) the primary vortex moves outward from the geometric axis. This displacement correlates with both the velocity trends (figure~\ref{fig:graph_velocity_of_R}b) and central curve movement derived from low-pressure regions (figures \ref{fig:3d_low_pressure_region}, \ref{fig:3d_central_line_projection}).

\begin{figure}[htbp]
	\centerline{\includegraphics[width=1\textwidth]{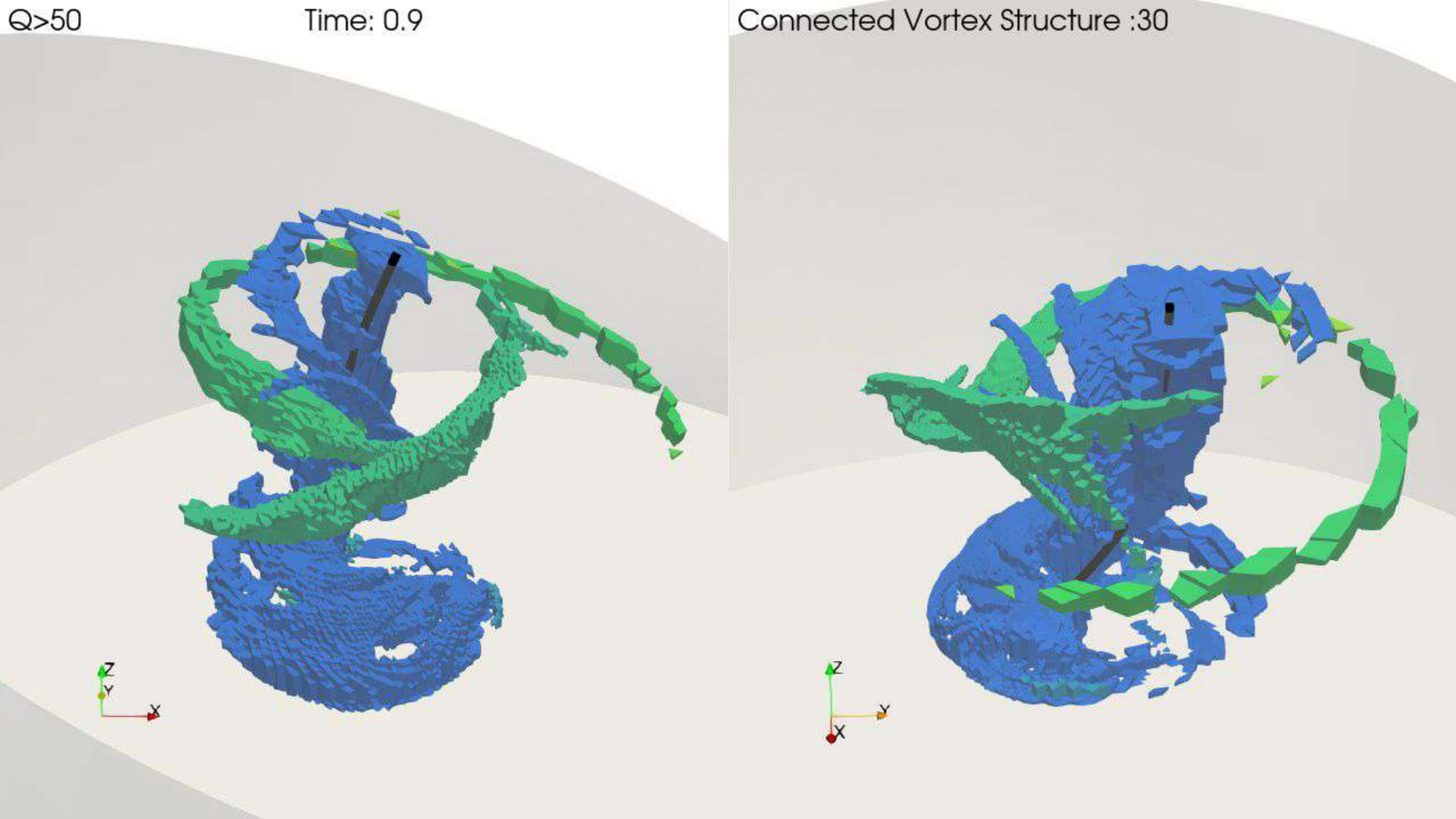}}
	\caption{3D vortex structure in curved domain $ \Omega_{C} $ $(R=1.5)$ at $t=0.9$. Front-(left) and side-views(right) showing Q-criterion vortices ($Q \geq 50$), with central curve (black). Colors distinguish individual vortex structures (connected components).}
\label{fig:3d_vortex_structure_Qcriterion_50}
\end{figure}

\begin{figure}[htbp]
	\centerline{\includegraphics[width=1\textwidth]{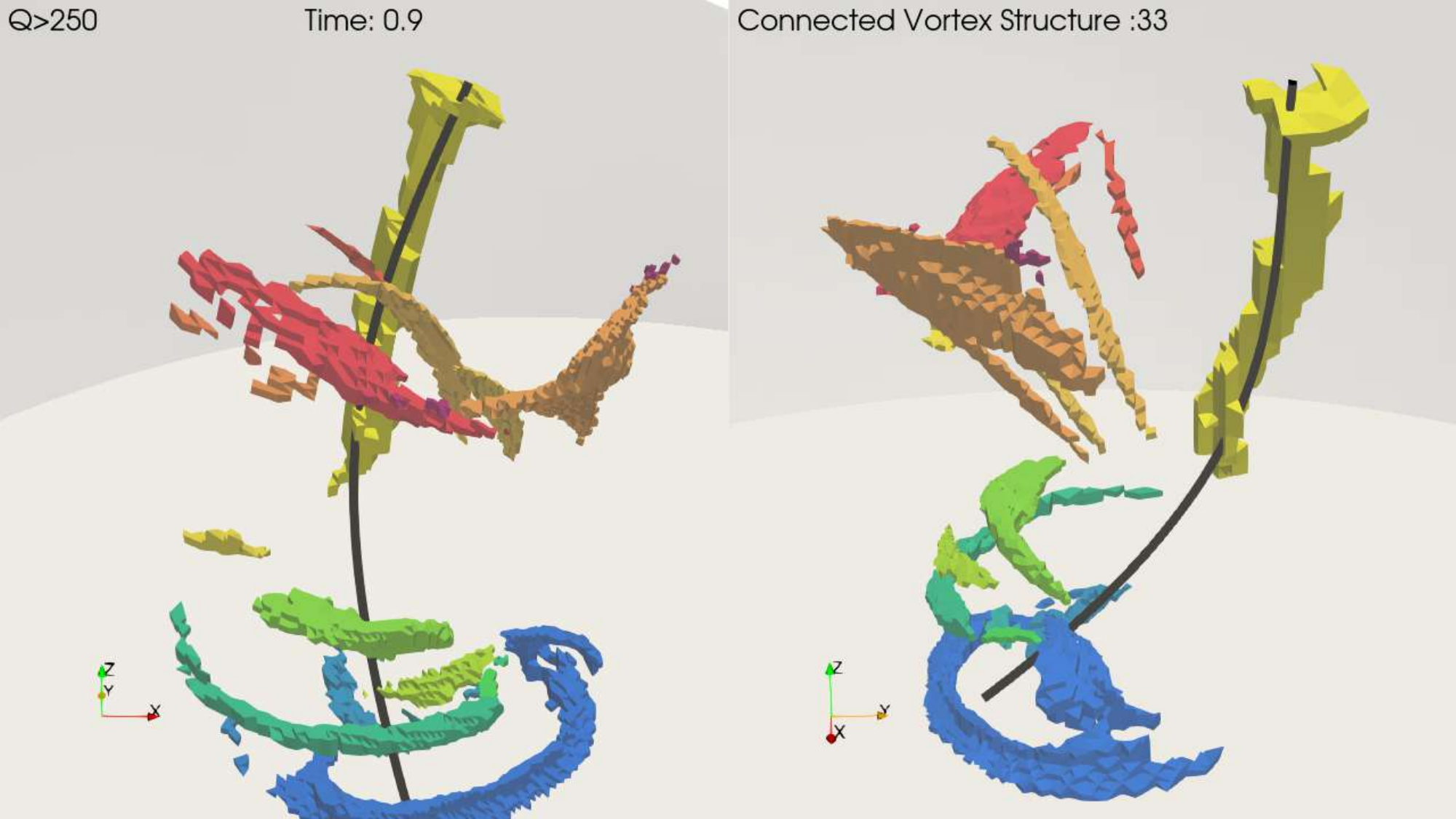}}
	\caption{3D vortex structure in curved domain $ \Omega_{C} $ $(R=1.5)$ at $t=0.9$. Front-(left) and side-views(right) showing Q-criterion vortices ($Q \geq 250$), with central curve (black). Colors distinguish individual vortex structures (connected components).}
\label{fig:3d_vortex_structure_Qcriterion_250}
\end{figure}

To better understand vortex dynamics, we visualize the connected structure of vortices (figures \ref{fig:3d_vortex_structure_Qcriterion_50}, \ref{fig:3d_vortex_structure_Qcriterion_250}). Each structure is shown in a unique color, with different structures represented by distinct colors. The figures also display the number of captured connected structures (called connected vortex structure). Analysis on vortex structures at lower thresholds ($ Q\geq50 $) shows that the primary vortex remains connected until $t=0.6$, and then splits at $t=0.9$ (figure~\ref{fig:3d_vortex_structure_Qcriterion_50}). At higher thresholds ($ Q\geq250 $), more complex vortex interactions are observed. At $t=0.9$ (figure~\ref{fig:3d_vortex_structure_Qcriterion_250}), two smaller vortices merge (shown in gold-brown), away from the primary vortex. This merged structures forms a stronger vortex ($ Q\geq750 $), captured in figure~\ref{fig:3d_vortex_Qcriterion_250_750}d, which satisfies our criteria for a secondary vortex.

Finally, we emphasize that \citet{BluesteinEtAl18} presents a radar study RaXPol on the El Reno, Oklahoma tornado (31 May 2013). Field observations shows that most long-lived secondary vortices developed within radius of maximum wind (RMW) and in the rear-left quadrant relative to the parent tornado. Our computational findings align closely with this real-world observations, lending credibility to the modeled mechanisms. In our simulations, secondary vortices develop in the rear-left region relative to the movement of the primary vortex and are located near the center of domain. The resemblance between the computational results and real tornado behavior strengthen the validity of our numerical simulation.

\section{Conclusion}\label{sec:conclusion}

Simulations in the curved cylindrical domain $ \Omega_{C} $ reveal dynamic behaviors that are absent in prior research by \citet{HsuNotsuYoneda16} in the straight cylindrical domain $ \Omega_{S} $. The maximum velocity location moves outward from the geometric axis, aligning with the displacement of the central curve derived from low-pressure regions. In contrast, in the straight cylindrical domain $ \Omega_{S} $, both the maximum velocity location and the central curve remain close to the geometric axis. Higher curvature accelerates this outward motion, supported by the angular momentum transfer across adjacent regions, where increased curvature leads to more complex and accelerated dynamics.

Also, vortices visualized using the Q-criterion highlight energy transfer between adjacent regions and clarify complex vortex dynamics in the curved domains. In particular, the numerical result in figure~\ref{fig:graph_energy_exchange} is qualitatively consistent with the picture of scale-local energy transfer \citep[Figure~1]{YonedaEtA2022} derived from turbulence snapshots in spatial averages at each time, that is, interaction of vortices between two adjacent scales. 

In the tornado study itself, the development of secondary vortex located left-rear quadrant relative to the primary vortex has also been documented in observational studies of El Reno tornado \citep{BluesteinEtAl18}, supporting the credibility of our numerical approach. The emergence of secondary vortex, identified at $ t=0.9 $ using a high Q-criterion threshold, occurs around the same time as the primary vortex moves outward from geometric axis. The movement of the central curve, derived from low-pressure regions further supports this observation. By observing velocity, pressure, angular momentum,
energy and the Q-criterion, these results enhance the reliability of vortex identification and improve understanding of vortex dynamics.

These findings underscore the significance of flow dynamics in curved geometries and highlight the potential to control flow properties through curvature. Future research could explore vertical behavior, external force, geometric modifications, or variations in physical properties in numerical simulations to study specific phenomena, 
especially, coherent structure of the energy cascade further. 

\textbf{Acknowledgements.}
Research of TY was partly supported by the JSPS Grants-in-Aid for Scientific Research JP24H00186.
HN was partly supported by the JSPS Grants-in-Aid for Scientific
Research JP20KK0058, JP20H01823, JP21H04431, JP24H00188 and JP25K00920, and JST CREST Grant Number JPMJCR2014.

\end{document}